\documentclass[reqno,a4paper,12pt]{amsart}
\usepackage[top=32mm, bottom=32mm, left=32mm, right=32mm]{geometry}
\usepackage{mathrsfs,enumerate}
\usepackage{amssymb}
\usepackage[colorlinks=true]{hyperref}

\newtheorem{theorem}{Theorem}

\theoremstyle{definition}
\newtheorem{definition}[theorem]{Definition}

\newtheorem{example}[theorem]{Example}

\newtheorem{question}[theorem]{Question}

\theoremstyle{remark}
\newtheorem{remark}[theorem]{Remark}

\begin{document}

\title[Sensitivity in general semi-flows]{Relationship among several types of sensitivity in general semi-flows}

\author[X. Wu]{Xinxing Wu}
\address[X. Wu]{School of Sciences, Southwest Petroleum University, Chengdu, Sichuan 610500, P.R. China}
\email{wuxinxing5201314@163.com}

\author[X. Zhang]{Xu Zhang}
\address[X. Zhang]{(Corresponding author)~Department of Mathematics, Shandong University, Weihai, Shandong 264209, P.R. China}
\email{xu$\_$zhang$\_$sdu@mail.sdu.edu.cn}

\subjclass[2010]{37B05, 54H20.}

\date{\today}

\keywords{Semi-flow, topological monoid, sensitivity.}

\begin{abstract}
In this paper, we show that there exists a monoid, on which neither the syndetic property nor the dual syndetic property holds, and there exists a strongly mixing  semi-flow with this monoid action which does not have thick sensitivity, syndetic sensitivity, thickly syndetic sensitivity, or thickly periodical sensitivity. Meanwhile, we show that there exists a thickly sensitive cascade which is not multi-sensitive. The first result answers positively Question 2,
and the first and the second results answer negatively Question 3 in \cite[A. Miller, A note about various types of sensitivity in general semiflows, Appl. Gen. Topol., 2018]{M2018-1}.
\end{abstract}

\maketitle

\section{Introduction}

In mathematics, a semigroup is an algebraic structure consisting of a set together with an associative binary operation.
The binary operation of a semigroup is most often denoted addition. A {\it monoid} is an algebraic structure intermediate
between groups and semigroups, and is a semigroup having an identity element, thus obeying all but one of the axioms of
a group; the existence of inverses is not required of a monoid.

Let $(X,d)$ be a metric space and $G$ be a non-compact abelian (commutative) topological monoid with the identity element $0$. A jointly continuous monoid action $\pi: G\times X\longrightarrow X$ of $G$ on the metric space $X$ is called a {\it semi-flow}
and denoted by $(G, X, \pi)$ or $(G, X)$. In particular, if the acting topological monoid is a topological group,
a semi-flow is called a {\it flow}. The element $\pi (t, x)$ will be denoted by $t. x$ or $tx$, so that the defining
conditions for a semi-flow have the following form:
$$
s.(t.x)=(s+t).x \text{ and } 0.x=x, \text{ for all } s, t\in G \text{ and all } x\in X.
$$
The maps
\begin{align*}
\pi_{t}: X &\longrightarrow X,\ t\in G\\
x &\longmapsto t.x
\end{align*}
are called {\it transition maps}. For any $x\in X$, the set $Gx:=\{t.x: t\in G\}$ is called the {\it orbit} of $x$.
A semi-flow $(G, X)$ is {\it minimal} if the orbit of every point $x\in X$ is dense in $X$, i.e., $\overline{Gx}=X$. Otherwise, it is called {\it non-minimal}.
For any subset $A$ of $X$ and any $t\in G$, let $t.A=\{t.x : x\in A\}$, denoted by $tA$ for convenience. 

In particular, let $\mathbb{N}=\left\{1, 2, 3, \ldots\right\}$ and $\mathbb{N}_0=\left\{0, 1, 2, \ldots\right\}$. Consider a continuous map $f: X\longrightarrow X$
and a monoid $G=\mathbb{N}_0$ (with the discrete topology), which leads to a natural semi-flow, that is, for any $x\in X$ and any $n\in \mathbb{N}_0$, $n. x=f^{n}(x)$.
This type of semi-flow is called a {\it cascade} and is often denoted by $(X, f)$ instead of $(\mathbb{N}_0, X)$.

Let $(X,d)$ be a metric space. A non-empty open subset of $X$ is said to be a {\it nopen}, or a {\it nopen} subset.
For any non-empty subset $U$ of $X$, the {\it diameter} $\mbox{diam}(U)$ of $U$ is
$\mbox{diam}(U)=\sup\{ d(x,y): x, y\in U\}$.

\begin{definition}\cite{M2018-1}
Let $G$ be a monoid. A subset $B$ of $G$ is
\begin{enumerate}[(1)]
\item {\it syndetic} if there exists a compact subset $K$ of $G$ such that, for every $t\in G$, $(t+K)\cap B\neq \emptyset$,
where $t+K=\left\{t+a: a\in K\right\}$; 
\item {\it thick} if, for every compact subset $K$ of $G$, there exists some $t\in G$ such that $t+K\subset B$;
\item {\it thickly syndetic} if, for every compact $K$ of $G$, there exists a syndetic subset $S$ of $G$ such that $S+K\subset B$;
\item {\it periodic} if there exist a closed syndetic sub-monoid $S$ of $G$ and $t\in G$ such that $t+S\subset B$;
\item {\it thickly periodic} if, for every compact subset $K$ of $G$, there exists a periodic subset $P$ of $G$ such that
$P+K\subset B$.
\end{enumerate}
\end{definition}

Clearly, a subset $B$ of $G$ is syndetic (thick) if and only if $G\setminus B$ is not thick (not syndetic).

\begin{definition} \cite[Definition 1.1]{M2018-1}
A monoid $G$ satisfies the {\it syndetic property} ({\it sp property}), or $G$ is an {\it sp monoid}, if all syndetic subsets of $G$
are non-compact; the {\it dual syndetic property} ({\it dsp property}), or a {\it dsp monoid} can be defined similarly: for every compact subset $K$ of $G$, the set $G\setminus K$ is a syndetic subset of $G$.
\end{definition}

\begin{definition} \cite[Definition 1.4]{M2018-1} A semi-flow $(G,X)$ is
\begin{enumerate}[(1)]
\item {\it strongly mixing} (StrM) if, for any two nopens $U,V$ of $X$, the set
$$D(U,V)=\{t\in G:~ tU\cap V\neq\emptyset\}$$
contains $G\setminus K$ for some compact subset $K$ of $G$;

\item {\it weak mixing} (WM) if, for any nopens $U_1,V_1,U_2,V_2$ of $X$, $$D(U_1,V_1)\cap D(U_2,V_2)\neq\emptyset;$$

\item {\it sensitive} (S) if there exists $\varepsilon>0$ such that for any nopen $U$ of $X$,
$$D(U,c)=\left\{t\in G:~ \mathrm{diam}(tU)>\varepsilon\right\}\neq\emptyset;$$

\item {\it strongly sensitive} (StrS) if there exists $\varepsilon>0$ such that for any nopen $U$ of $X$,
$D(U,c)$ contains $G\setminus K$ for some compact subset $K$ of $G$;

\item {\it multi-sensitive} (MulS) if there exists $\varepsilon>0$ such that for any $n\in \mathbb{N}$ and
any nopens $U_1, \ldots, U_n$ of $X$, $\cap_{i=1}^{n}D(U_i, \varepsilon) \neq \emptyset$;

\item {\it strongly multi-sensitive} (StrMulS) if there exists $\varepsilon>0$ such that for any $n\in \mathbb{N}$
and any nopens $U_1, \ldots,U_n$ of $X$, $\cap_{i=1}^{n}D(U_i, \varepsilon)$ contains $G\setminus K$ for some compact subsets $K$ of $G$;

\item {\it thickly sensitive} (TS) if there exists $\varepsilon>0$ such that for any nopen $U$ of $X$,
$D(U,\varepsilon)$ is a thick subset of $G$;

\item {\it syndetically sensitive} (SyndS) if there exists $\varepsilon>0$ such that for any nopen $U$ of $X$,
$D(U,\varepsilon)$ is a syndetic subset of $G$;

\item {\it thickly syndetically sensitive} (TSynS) if there exists $\varepsilon>0$ such that for any nopen $U$ of $X$,
$D(U, \varepsilon)$ is a thickly syndetic subset of $G$;

\item {\it periodically sensitive} (PerS) if there exists $\varepsilon>0$ such that for any nopen $U$ of $X$,
$D(U,\varepsilon)$ is a periodic subset of $G$;

\item {\it thickly periodically sensitive} (TPerS) if there exists $\varepsilon>0$ such that for any nopen $U$ of $X$,
$D(U,\varepsilon)$ is a thickly periodic subset of $G$.
\end{enumerate}
\end{definition}

Ceccherini-Silberstein and Coornaert proved that every transitive semi-flow on an infinite Hausdorff uniform space
admitting a dense set of periodic points is sensitive~\cite{CC2013}. Miller obtained that every non-minimal syndetically transitive semi-flow is syndetically sensitive~\cite{M2015}, generalizing the main results in \cite{BBCDS1992,M2007,WLF2012}. Further, Miller gave a summary on the sensitivity of semi-flows
with monoid actions \cite{M2017}. Money conducted a systemic investigation on chaotic properties for semi-flows~\cite{Mo2015}. Wang showed that every
$M$-action on a compact Hausdorff uniform space having at least two disjoint compact invariant subsets is thickly syndetically sensitive~\cite{W2018}.
Recently, Miller studied the relationships among various types of sensitivity in general
semi-flows and proposed several open questions~\cite{M2018-1}. For more recent results on sensitivity, refer to \cite{IP2019,KM2008,LOW2017,M2011,M2018-2,WC2017,WLML2019,WOC2016,WWC2015,WZ2012}
and some references therein.

\begin{question}\label{Q-2}{\rm \cite[Question~2]{M2018-1}
Find examples showing that in the implications StrS $\Rightarrow$ TS
and StrS $\Rightarrow$ TSyndS the condition (sp) is indeed needed,
and that in the implication SM $\Rightarrow$ SyndS the condition (dsp) is indeed needed.}
\end{question}

\begin{question}\label{Q-3}{\rm \cite[Question~3]{M2018-1}
Investigate if, in general, StrS implies TPerS, TS implies MulS, and TS implies SyndS.}
\end{question}

\begin{remark}\label{R-1}
Recently, we \cite{WMCL} proved that
\begin{enumerate}[(1)]
\item there exist two non-syndetically sensitive cascades defined on complete metric spaces whose product is cofinitely sensitive;

\item there exists a syndetically sensitive semi-flow $(G, X)$ defined on a complete metric space $X$ such that
$(G_1, X)$ is not sensitive for some syndetic closed sub-monoid $G_1$ of $G$, which gives a negative answer to \cite[Question~43]{M2017};

\item there exists a thickly sensitive cascade which is not syndetically sensitive, i.e., (TS) $\nRightarrow$ (SyndS).
\end{enumerate}
\end{remark}

This paper proves that there exists a monoid, on which neither the syndetic property nor the dual syndetic property holds, and there exists a semi-flow
given by this monoid satisfying that the strong mixing property can not imply the thick sensitivity, the syndetic sensitivity, the thickly syndetic sensitivity,
or the thickly periodical sensitivity (See Example~\ref{Exa-1}). Meanwhile, it is proved that there exists a thickly sensitive cascade which is not multi-sensitive
(See Example~\ref{Exa-2}). These results, together with (3) of Remark~\ref{R-1}, give complete answers to Questions \ref{Q-2} and \ref{Q-3}.

\section{Relationship among several types of sensitivity}

\subsection{A monoid without the sp or dsp property}

In this subsection, a monoid is provided, on which neither the syndetic property nor the dual syndetic property holds, and there exists a semi-flow
given by this monoid satisfying that the strong mixing property can not imply the thick sensitivity, the syndetic sensitivity, the thickly syndetic sensitivity,
or the thickly periodical sensitivity.

\begin{example}\label{Exa-1}
Consider the set $G=\mathbb{N}_0 \cup \{\infty\}$ with the discrete topology, which can be thought of as the one-point compactification of $\mathbb{N}_0$.
The addition `$+$' defined on $G$ is given as follows:
$$
x+y=
\begin{cases}
x+y, & \text{if } x,y\in\mathbb{N}_0, \\
\infty, & \text{if } x=\infty \text{ or } y=\infty.
\end{cases}
$$

First, we verify the following properties.
\begin{itemize}
\medskip
\item[(1)] It is clear that $G$ is a monoid with the identity element $0$.

\medskip

\item[(2)] $G$ does not satisfy the dsp or the sp property.

\medskip

Take a compact subset $K=\{\infty\}$ of $G$. It suffices to show that $G\setminus K$ is not syndetic,
implying that $G$ does not satisfy the dsp property. It can be verified that for any non-empty compact subset $K_1$ of $G$,
$\infty+K_1=\{\infty\}$. Then,
$$
(\infty+K_1)\cap(G\setminus K)=\{\infty\}\cap (G\setminus\{\infty\})=\emptyset.
$$
Thus, $G\setminus K$ is not syndetic.

This, together with \cite[Proposition~1.3]{M2018-1}, implies that $G$ does not satisfy the sp property.

\medskip

\item [(2$^{'}$)] Every thick subset of $G$ contains $\infty$.

\medskip

Take any thick subset $K$ of $G$, from the definition, it follows that for the compact subset $\{\infty\}$, there exists
$g\in G$ such that $g+\{\infty\}=\{\infty\}\subset K$.

\medskip

\item [(2$^{''}$)] Every thickly periodical subset of $G$ contains $\infty$.

\medskip

Take any thickly periodical subset $K$, for the compact subset $\{\infty\}$,
there exists a periodic subset $P$ of $G$ such that $P+\{\infty\}=\{\infty\}\subset K$.

\medskip

\item[(2$^{'''}$)] Every thickly syndetic subset of $G$ contains $\infty$.

\medskip

Take any thickly syndetic subset $K$, for the compact subset $\{\infty\}$,
there exists a syndetic subset $S\subset G$ such that
$S+\{\infty\}=\{\infty\}\subset K$.
\end{itemize}

Now, we construct a semi-flow. Let $I=[0,1]$ and $f(x)=1-|1-2x|$ be the tent map from $I$ to $I$.
It is easy to see that $(I, f)$ is strongly mixing. Define the map $\pi:G\times X\longrightarrow X$ as follows:
\begin{itemize}
\item[(3)] $\pi(n,x)=f^n(x)$, $\forall n\in\mathbb{N}_0$ and $\forall x\in I$;
\item[(4)] $\pi(\infty,x)=0$, $\forall x\in I$.
\end{itemize}
It can be verified that
\begin{itemize}
\item[(5)] $(G,I,\pi)$ is a semi-flow;
\item[(6)] $(G,I,\pi)$ is strongly mixing.
\end{itemize}

For any nopen subsets $U,V$ of $I$, since $(I,f)$ is strongly mixing, noting that
$$
G\setminus\{g\in G:\ gU\cap V\neq\emptyset\} \subset \{\infty\}\cup\left(\mathbb{N}_0\setminus\{n\in\mathbb{N}_0:\ f^n(U)\cap V\neq\emptyset\}\right),
$$
we have that $G\setminus\{g\in G:\ gU\cap V\neq\emptyset\}$ is a compact subset of $G$, that is, $(G,I,\pi)$ is strongly mixing.

\medskip

Next, we verify the following properties.

\begin{itemize}

\medskip

\item[(7)] From (6) and \cite[Proposition 2.1]{M2018-1}, it follows that $(G,I,\pi)$ is strongly sensitive.

\medskip

\item[(8)] $(G,I,\pi)$ is not thickly sensitive.

\medskip

For any $0<\varepsilon<1$, it is easy to see that $\{g\in G:\ \mbox{diam}(gI)<\varepsilon\}=\{\infty\}$. This, together with (2$^{'}$),
implies that $\{g\in G:\ \mbox{diam}(gI)\geq \varepsilon\}=G\setminus \{\infty\}$ is not thick. Thus, $(G,I,\pi)$ is not thickly sensitive.

\medskip

\item[(9)]$(G,I,\pi)$ is not syndetically sensitive.

\medskip

Applying similar arguments as in (8), and the fact that $G\setminus\{\infty\}$
is not syndetic in (2) yields that $(G,I,\pi)$ is not syndetically sensitive.

\medskip

\item[(10)]$(G,I,\pi)$ is not thickly periodically sensitive.

\medskip

Applying similar discussions as in (8), this together with (2$^{''}$), yields that
$(G,I,\pi)$ is not thickly periodically sensitive.

\medskip

\item[(11)]$(G,I,\pi)$ is not thickly syndetically sensitive.

\medskip

Applying similar discussions as in (8), this together with (2$^{'''}$), yields that
$(G,I,\pi)$ is not thickly syndetically sensitive.
\end{itemize}
\end{example}

\begin{remark}
\begin{enumerate}[(1)]
\item Example~\ref{Exa-1} gives a positive answer to Question~\ref{Q-2}.
\item Miller \cite{M2019} obtained a weakly mixing semi-flow defined on a non-compact metric space which is not thickly sensitive.
Example~\ref{Exa-1} shows that there exists a strongly mixing semi-flow defined on a compact metric space which does not
have thick sensitivity, syndetic sensitivity, thickly syndetic sensitivity, or thickly periodical sensitivity.

\item Clearly, $G_1=\{0, \infty\}$ is a syndetic closed sub-monoid of $G$ in Example~\ref{Exa-1}. Meanwhile, it can
be verified that $(G_1, I, \pi)$ is not sensitive. Similarly to (2) of Remark~\ref{R-1}, this, together with the strong sensitivity of $(G, I, \pi)$, also
gives a negative answer to \cite[Question~43]{M2017}.
\end{enumerate}
\end{remark}

\subsection{A thickly sensitive, but not multi-sensitive, semi-flow}

Similarly to the construction of \cite[Example~4]{WMCL},
a thickly sensitive semi-flow which is not multi-sensitive is constructed in this subsection.

\begin{example}\label{Exa-2}
Let
$$L_0=\mathscr{L}_0=0,\  L_1=\mathscr{L}_1=2,$$
and
$$L_n=2^{L_1+ \cdots +L_{n-1}}\cdot (2n),\ \mathscr{L}_{n}=L_1+L_2+ \cdots +L_n,\ n\geq2,$$
and take
$$
\begin{aligned}
X=\left[0.5, 1.5\right]&\bigcup \left(\bigcup_{n=1}^{+\infty}\bigcup_{i=0}^{2^{\mathscr{L}_{2n-2}}}\left[\mathscr{L}_{2n-1}+i,
\mathscr{L}_{2n-1}+i+\tfrac{1}{2n}\right]\right)\\
&\bigcup \left(\bigcup_{n=1}^{+\infty}\bigcup_{i=0}^{2^{\mathscr{L}_{2n-1}}}
\left[\mathscr{L}_{2n}+(2n+1)i, \mathscr{L}_{2n}+(2n+1)i+2n\right]\right),
\end{aligned}
$$
and
$$
\begin{aligned}
Y=\left[0.5, 1.5\right]&\bigcup \left(\bigcup_{n=1}^{+\infty}\bigcup_{i=0}^{2^{\mathscr{L}_{2n-2}}}
\left[\mathscr{L}_{2n-1}+2ni,
\mathscr{L}_{2n-1}+2ni+2n-1\right]\right)\\
&\bigcup \left(\bigcup_{n=1}^{+\infty}\bigcup_{i=0}^{2^{\mathscr{L}_{2n-1}}}
\left[\mathscr{L}_{2n}+i, \mathscr{L}_{2n}+i+\tfrac{1}{2n+1}\right]\right).
\end{aligned}
$$
For $n\in \mathbb{N}$, let
$$
\begin{aligned}
\mathscr{A}_n&=\left[\mathscr{L}_{2n-1}+2^{\mathscr{L}_{2n-2}},
\mathscr{L}_{2n-1}+2^{\mathscr{L}_{2n-2}}+\tfrac{1}{2n-1}\right],\\
\mathscr{B}_n&=\left[\mathscr{L}_{2n}+(2n+1)\cdot 2^{\mathscr{L}_{2n-1}},
\mathscr{L}_{2n}+(2n+1)\cdot 2^{\mathscr{L}_{2n-1}}+2n\right],\\
\mathscr{C}_n&=\left[\mathscr{L}_{2n-1}+2n\cdot 2^{\mathscr{L}_{2n-2}},
\mathscr{L}_{2n-1}+2n\cdot 2^{\mathscr{L}_{2n-2}}+(2n-1)\right],\\
\mathscr{D}_{n}&=\left[\mathscr{L}_{2n}+2^{\mathscr{L}_{2n-1}},
\mathscr{L}_{2n}+2^{\mathscr{L}_{2n-1}}+\tfrac{1}{2n+1}\right].
\end{aligned}
$$
Define two linear maps $f: X\to X$ and $g: Y\to Y$ defined by
$$
f(x)=\left\{
\begin{aligned}
&0.5x+1.75, \text{ if }  x\in [0.5, 1.5], \\
&x+1, \text{ if }  x\in \left[\mathscr{L}_{2n-1}+i, \mathscr{L}_{2n-1}+i+\tfrac{1}{2n}\right]
\text{ for some } 0\leq i<2^{\mathscr{L}_{2n-2}}, n\in \mathbb{N}, \\
&x+2n+1, \text{ if }  x\in \left[\mathscr{L}_{2n}+(2n+1)i, \mathscr{L}_{2n}+(2n+1)i+2n\right]\\
&\quad\quad\quad\quad\quad\quad \text{ for some } 0\leq i<2^{\mathscr{L}_{2n-1}}, n\in \mathbb{N}, \\
&2n(2n-1)\left(x-
\mathscr{L}_{2n-1}-2^{\mathscr{L}_{2n-2}}\right)+\mathscr{L}_{2n}, \text{ if } x\in \mathscr{A}_n, n\in \mathbb{N},\\
&\tfrac{1}{2n(2n+1)}\left(x-\mathscr{L}_{2n}-(2n+1)\cdot 2^{\mathscr{L}_{2n-1}}\right)+\mathscr{L}_{2n+1},
\text{ if } x\in \mathscr{B}_n, n\in \mathbb{N},
\end{aligned}
\right.
$$
and
$$
g(x)=\left\{
\begin{aligned}
&x+1.5, \text{ if }  x\in [0.5, 1.5], \\
&x+1, \text{ if }  x\in \left[\mathscr{L}_{2n}+i, \mathscr{L}_{2n}+i+\tfrac{1}{2n}\right]
\text{ for some } 0\leq i<2^{\mathscr{L}_{2n-1}}, n\in \mathbb{N}, \\
&x+2n, \text{ if }  x\in \left[\mathscr{L}_{2n-1}+2ni, \mathscr{L}_{2n-1}+2ni+(2n-1)\right]\\
&\quad\quad\quad\quad\text{ for some } 0\leq i<2^{\mathscr{L}_{2n-2}}, n\in \mathbb{N}, \\
&\tfrac{1}{(2n+1)\cdot (2n-1)}\left(x-\mathscr{L}_{2n-1}-2n\cdot 2^{\mathscr{L}_{2n-2}}\right)+\mathscr{L}_{2n},
\text{ if } x\in \mathscr{C}_n, n\in \mathbb{N},\\
&(2n+1)^{2} \left(x-\mathscr{L}_{2n}-2^{\mathscr{L}_{2n-1}}\right)+\mathscr{L}_{2n+1},
\text{ if } x\in \mathscr{D}_n, n\in \mathbb{N},
\end{aligned}
\right.
$$
respectively. Clearly, $f$ and $g$ are continuous. Rearrange all closed intervals as follows:
$$[0.5, 1.5], [\mathscr{L}_1, \mathscr{L}_1+\tfrac{1}{2}], \ldots,
[\mathscr{L}_{1}+2^{\mathscr{L}_0}, \mathscr{L}_{1}+2^{\mathscr{L}_0}+\tfrac{1}{2}],$$ $$[\mathscr{L}_2, \mathscr{L}_{2}+2], \ldots, [\mathscr{L}_2+2\cdot 2^{\mathscr{L}_1},
\mathscr{L}_2+2\cdot 2^{\mathscr{L}_1}+2],$$ $$\vdots$$
$$[\mathscr{L}_{2n-1}, \mathscr{L}_{2n-1}+\tfrac{1}{2n}], \ldots, [\mathscr{L}_{2n-1}+2^{\mathscr{L}_{2n-2}},
\mathscr{L}_{2n-1}+2^{\mathscr{L}_{2n-2}}+\tfrac{1}{2n}],$$ $$[\mathscr{L}_{2n}, \mathscr{L}_{2n}+2n], \ldots, [\mathscr{L}_{2n}+(2n+1)2^{\mathscr{L}_{2n-1}},
\mathscr{L}_{2n}+(2n+1)2^{\mathscr{L}_{2n-1}}+2n], \ldots$$ of $X$ by this natural order and denote by  $$I_0, I_1, I_2, \ldots.$$
It is easy to see that
$$I_{2k+2+2^{\mathscr{L}_1}+\cdots +2^{\mathscr{L}_{2k-1}}}=[\mathscr{L}_{2k+1}, \mathscr{L}_{2k+1}+\tfrac{1}{2(k+1)}],$$
$$\vdots$$
$$I_{2k+2+2^{\mathscr{L}_1}+\cdots +2^{\mathscr{L}_{2k}}}=[\mathscr{L}_{2k+1}+2^{\mathscr{L}_{2k}},
\mathscr{L}_{2k+1}+2^{\mathscr{L}_{2k}}+\tfrac{1}{2(k+1)}],$$ $$I_{2k+3+2^{\mathscr{L}_1}+\cdots +2^{\mathscr{L}_{2k}}}=
[\mathscr{L}_{2k+2}, \mathscr{L}_{2k+2}+2(k+1)],$$
$$\vdots$$
$$I_{2k+3+2^{\mathscr{L}_1}+\cdots +2^{\mathscr{L}_{2k+1}}}=
[\mathscr{L}_{2k+2}+(2k+3)\cdot 2^{\mathscr{L}_{2k+1}}, \mathscr{L}_{2k+2}+(2k+3)\cdot 2^{\mathscr{L}_{2k+1}}+2(k+1)].$$

Similarly, rearrange all closed intervals of $Y$ by this natural order and denote by
$$J_0, J_1, J_2, \ldots.$$ It is easy to see that
$$J_{2k+2+2^{\mathscr{L}_1}+\cdots +2^{\mathscr{L}_{2k-1}}}=[\mathscr{L}_{2k+1}, \mathscr{L}_{2k+1}+2k+1],$$
$$\vdots$$
$$J_{2k+2+2^{\mathscr{L}_1}+\cdots +2^{\mathscr{L}_{2k}}}=[\mathscr{L}_{2k+1}+2\cdot (k+1)\cdot 2^{\mathscr{L}_{2k}},
\mathscr{L}_{2k+1}+2\cdot (k+1)\cdot 2^{\mathscr{L}_{2k}}+2k+1],$$ $$J_{2k+3+2^{\mathscr{L}_1}+\cdots +2^{\mathscr{L}_{2k}}}=
[\mathscr{L}_{2k+2}, \mathscr{L}_{2k+2}+\tfrac{1}{2k+3}],$$
$$\vdots$$ $$J_{2k+3+2^{\mathscr{L}_1}+\cdots +2^{\mathscr{L}_{2k+1}}}=
[\mathscr{L}_{2k+2}+2^{\mathscr{L}_{2k+1}}, \mathscr{L}_{2k+2}+2^{\mathscr{L}_{2k+1}}+\tfrac{1}{2k+3}].$$
Note that $f$ is a linear homeomorphism from $I_{n}$ to $I_{n+1}$, and $g$ is also a linear homeomorphism from $J_{n}$ to $J_{n+1}$,
$n\in \mathbb{N}_0$. According to the construction of $f$ and $g$, it can be verified that
\begin{enumerate}[(i)]

\item\label{i} $f$ and $g$ are continuous;

\item\label{iv} for any $n\in \left[2^{\mathscr{L}_1}+\cdots + 2^{\mathscr{L}_{2k-1}}+2k+2, 2^{\mathscr{L}_1}+\cdots +
2^{\mathscr{L}_{2k}}+2k+2\right]$ ($k\in \mathbb{N}$),
$$
f^{n}([0.5, 1.5])=\left[\mathscr{L}_{2k+1}+n-a_k, \mathscr{L}_{2k+1}+n-a_k+\tfrac{1}{2k+2}\right],
$$
and
$$
g^n([0.5, 1.5])=\left[\mathscr{L}_{2k+1}+2(k+1)(n-a_k), \mathscr{L}_{2k+1}+2(k+1)(n-a_k)+2k+1\right],
$$
where $a_k=2^{\mathscr{L}_1}+\cdots + 2^{\mathscr{L}_{2k-1}}+2k+2$,
implying that $$\mathrm{diam}(f^{n}([0.5,1.5]))=\tfrac{1}{2k+2}\ \mbox{and}\ \mathrm{diam}(g^{n}([0.5,1.5]))=2k+1;$$

\item\label{v} for any $n\in \left[2^{\mathscr{L}_1}+\cdots + 2^{\mathscr{L}_{2k}}+2k+3, 2^{\mathscr{L}_1}+\cdots +
2^{\mathscr{L}_{2k+1}}+2k+3\right]$ ($k\in \mathbb{N}$),
$$
f^{n}([0.5, 1.5])=\left[\mathscr{L}_{2k+2}+(n-b_k)(2k+3), \mathscr{L}_{2k+2}+(n-b_k)(2k+3)+2k+2\right],
$$
and
$$
g^{n}([0.5, 1.5])=\left[\mathscr{L}_{2k+2}+n-b_k, \mathscr{L}_{2k+2}+n-b_k+\tfrac{1}{2k+3}\right],
$$
where $b_k=2^{\mathscr{L}_1}+\cdots +
2^{\mathscr{L}_{2k}}+2k+3$, implying that $$\mathrm{diam}(f^{n}([0,1]))=2k+2\ \mbox{and}\
\mathrm{diam}(g^{n}([0.5,1.5]))=\tfrac{1}{2k+3}.$$
\end{enumerate}
Take $Z=\left(X\times \{0\}\right) \cup \left(\{0\}\times Y\right)\subset \mathbb{R}^2$. It is easy to see that
$Z$ is a complete metric subspaces of $\mathbb{R}^2$.
Define $F: Z\longrightarrow Z$ by
$$
F(x, y)=
\begin{cases}
(f(x), 0), &\text{if } (x, y)\in X\times \{0\}, \\
(0, g(y)), &\text{if } (x, y)\in \{0\}\times Y.
\end{cases}
$$
Clearly, $F$ is continuous as $f$ and $g$ are continuous.

\medskip

{\bf Claim 1.} $(Z, F)$ is thickly sensitive.

\medskip

Given any nopen subset $W$ of $Z$, there exists a nopen subset $U \subset X$
or $V\subset Y$ such that $U\times \{0\}\subset W$ or $\{0\}\times V\subset W$. Without loss of generality,
assume $U\times \{0\}\subset W$. It follows from \eqref{iv}--\eqref{v} that
there exist non-degenerate closed interval $[\alpha, \beta] \subset [0.5, 1.5]$
and $p\in \mathbb{N}_0$ such that $f^p([\alpha, \beta])\subset U$.
Recall that $f$ is a piecewise linear mapping, for any $n \in \mathbb{N}_0$, one has
$$
\mathrm{diam}(f^n([\alpha, \beta]))=(\beta-\alpha)\cdot \mathrm{diam}(f^n([0.5, 1.5])),
$$
implying that
$$
\begin{aligned}
\mathrm{diam}(F^n(W))&\geq \mathrm{diam}(F^n(U\times \{0\}))
\geq \mathrm{diam}(F^n(f^p([\alpha, \beta])\times \{0\}))\\
&=\mathrm{diam}(f^{n+p}([\alpha, \beta]))
=(\beta-\alpha)\cdot \mathrm{diam}(f^{n+p}([0.5, 1.5])).
\end{aligned}
$$
This, together with \eqref{v}, implies that
there exists $K\in \mathbb{N}$ such that
$$
\bigcup_{k=K}^{+\infty}\left[2^{\mathscr{L}_1}+\cdots + 2^{\mathscr{L}_{2k}}+2k+3-p, 2^{\mathscr{L}_1}+\cdots +
2^{\mathscr{L}_{2k+1}}+2k+3-p\right]\subset \mathrm{D}(U, 1),
$$
i.e., $\mathrm{D}(U, 1)$ is a thick set. Therefore, $F$ is thickly sensitive.

\medskip

{\bf Claim 2.} $(Z, F)$ is not multi-sensitive.

\medskip

For any $0<\varepsilon<1$, take two open subsets $U_1=[0.5, 0.5+\tfrac{\varepsilon}{2})\times \{0\}$ and
$U_2=\{0\}\times [0.5, 0.5+\tfrac{\varepsilon}{2})$ of $Z$. From the proof of Claim 1 and \eqref{iv}--\eqref{v},
it follows that
\begin{enumerate}[(1)]
\item for any $n\in \left[2^{\mathscr{L}_1}+\cdots + 2^{\mathscr{L}_{2k-1}}+2k+2, 2^{\mathscr{L}_1}+\cdots +
2^{\mathscr{L}_{2k}}+2k+2\right]$,
$$
\mathrm{diam}(F^{n}(U_1))=\mathrm{diam}\left(f^n\left(\left[0.5, 0.5+\tfrac{\varepsilon}{2}\right)\right)\right)
=\tfrac{\varepsilon}{4(k+1)}<\varepsilon;
$$
\item for any $n\in \left[2^{\mathscr{L}_1}+\cdots + 2^{\mathscr{L}_{2k}}+2k+3, 2^{\mathscr{L}_1}+\cdots +
2^{\mathscr{L}_{2k+1}}+2k+3\right]$,
$$
\mathrm{diam}(F^{n}(U_2))=\mathrm{diam}\left(g^n\left(\left[0.5, 0.5+\tfrac{\varepsilon}{2}\right)\right)\right)
=\tfrac{\varepsilon}{2(2k+3)}<\varepsilon.
$$
\end{enumerate}
This implies that
$$
D(U_1, \varepsilon) \cap D(U_2, \varepsilon)=\emptyset.
$$
Thus, $(Z, F)$ is not multi-sensitive.
\end{example}

\begin{remark}
From Examples~\ref{Exa-1} and \ref{Exa-2}, it follows that (StrS) $\nRightarrow$ (TPerS) and (TS) $\nRightarrow$ (MulS).
This, together with (3) of Remark~\ref{R-1}, answers negatively Question~\ref{Q-3}.
\end{remark}

\section*{Acknowledgments}
{\footnotesize This work was supported by the National Natural Science Foundation of China (Nos. 11601449 and 11701328),
the National Natural Science Foundation of China (Key Program) (No. 51534006), Science and Technology Innovation Team of Education
Department of Sichuan for Dynamical System and its Applications (No. 18TD0013), Youth Science and Technology Innovation Team of
Southwest Petroleum University for Nonlinear Systems (No. 2017CXTD02), Scientific Research Starting Project of Southwest Petroleum
University (No. 2015QHZ029), Shandong Provincial Natural Science Foundation, China (Grant ZR2017QA006), and Young Scholars Program
of Shandong University, Weihai (No. 2017WHWLJH09).}

\baselineskip=2pt

\providecommand{\bysame}{\leavevmode\hbox to3em{\hrulefill}\thinspace}
\providecommand{\MR}{\relax\ifhmode\unskip\space\fi MR }
\providecommand{\MRhref}[2]{%
  \href{http://www.ams.org/mathscinet-getitem?mr=#1}{#2}
}
\providecommand{\href}[2]{#2}

\end{document}